\documentclass{amsart}

\usepackage{amsmath}
\usepackage{amssymb}
\usepackage{faktor}
\usepackage{graphicx}
\usepackage{tikz-cd}
\usepackage{enumitem}
\usepackage{tikz}
\usetikzlibrary{intersections}

\newcommand{\N}{\ensuremath{\mathbb{N}}}   % natural numbers 
\newcommand{\C}{\ensuremath{\mathbb{C}}}   % complex numbers

\renewcommand{\epsilon}{\varepsilon}       % normal epsilon and
\renewcommand{\phi}{\varphi}		        % phi
\renewcommand{\k}{\ensuremath{\mathcal{K}}}  %compact operators

\newcommand{\T}{\ensuremath{\mathcal{T}}} %Toeplitz algebra

\DeclareMathOperator{\id}{id} %identity

\newtheorem*{thm}{Theorem}
\newtheorem*{lem}{Lemma}

\theoremstyle{definition}
\newtheorem*{dfn}{Definition}
\newtheorem*{ex}{Example}
\theoremstyle{remark}

\author{Laura Brake}
\email{l.brake@uni-muenster.de} 

\author{Wilhelm Winter}
\email{wwinter@uni-muenster.de} 
\address{Mathematisches Institut der WWU M\"unster, Einsteinstr. 62, 48149 M\"unster, Germany} 

\thanks{Supported by the Deutsche Forschungsgemeinschaft through SFB 878}

\title{The Toeplitz algebra has nuclear dimension one}

\begin{document}

\begin{abstract}
We prove the title by constructing 2-colourable completely positive approximations for the Toeplitz algebra. Besides results about nuclear dimension and completely positive contractive order zero maps, our argument involves projectivity of the cone over a finite dimensional  $\textup{C}^*$-algebra and Lin's theorem on almost normal matrices. 
\end{abstract}

\maketitle

Nuclear dimension is a notion of covering dimension for amenable C$^*$-algebras; it was introduced in \cite{WinZac:dimnuc} and, in contrast to the decomposition rank of \cite{KirchbergWinter:CovDim}, may take finite values also for infinite C$^*$-algebras. The most basic example of an infinite C$^*$-algebra is the Toeplitz algebra $\mathcal{T}$, which may be thought of as the C$^*$-algebra generated by the unilateral shift on $\ell^2(\mathbb{N})$; it is an extension of the continuous functions on the circle by the compact operators on $\ell^2(\mathbb{N})$:
\[
0 \longrightarrow \mathcal{K} \longrightarrow \mathcal{T} \longrightarrow C(S^1) \longrightarrow 0
\]
Finiteness of nuclear dimension is preserved under taking extensions, and it was observed already in \cite[Example~6.3]{WinZac:dimnuc} that the Toeplitz algebra has nuclear dimension either 1 or 2. The exact value remained unknown; cf.\ \cite[Problem~9.2]{WinZac:dimnuc}. 

Before settling the matter, let us explain why we even care. The arguably most striking appearance of nuclear dimension was in Elliott's classification programme for simple nuclear C$^*$-algebras; see  \cite{Win:ICM} for an overview. In this context, there now is a result as complete as can be, and the crucial hypothesis is just finiteness of nuclear dimension, with no reference to the actual value. In hindsight the latter is not too surprising, since we now know that for separable, simple, unital C$^*$-algebras the only possible values for the nuclear dimension are $0$, $1$, and infinity. On the other hand, in the nonsimple case the nuclear dimension potentially carries more information (after all, for continuous trace C$^*$-algebras it coincides on the nose with covering dimension of the spectrum). 
But even then one could argue why the exact value is relevant for an example which is otherwise as well understood as the Toeplitz algebra. Indeed, we do not so much expect new insights in the structure of said algebra, but our construction is a new instance of \emph{dimension reduction phenomena} in noncommutative topology.

Dimension reduction occurs in its most simple form in the commutative setting, where dimension is additive, and not multiplicative, under taking products: for example, an interval may be approximated by open covers with two colours (where members of the cover with the same colour are not allowed to intersect), and hence has dimension at most one. (Showing that it is not zero, i.e.\ finding a lower bound, is another matter.) But then the square, being a product of two intervals, can obviously be approximated by open covers with \emph{four} colours, hence has dimension at most three. To improve this and find approximations with \emph{three} colours (and to arrive at dimension two, as it should be) one has to employ (an ever so small amount of) geometric insight and come up, for instance, with a bricklaying scheme. 
 
Nuclear dimension expresses all this in terms of functions, with open covers replaced by partitions of unity and Cartesian products by tensor products. For commutative C$^*$-algebras nothing changes, but the noncommutative situation becomes more complicated. On the one hand, we do not know in general whether nuclear dimension is subadditive with respect to tensor products. On the other hand, there are several situations where dimension reduces when there is enough `noncommutative space' available, in particular \cite{Gong:reduction, KirchbergRordam:pi3, TikWin:Z-dr, MS:DMJ, SWW:Invent, BBSTWW:arXiv}. In all of these, dimension reduction is quite dramatic, and eventually relies (more or less explicitly) on the existence of simple -- in most cases, \emph{strongly self-absorbing} -- local tensor factors, like the Cuntz algebras $\mathcal{O}_2$ and $\mathcal{O}_\infty$, UHF algebras, or the Jiang--Su algebra $\mathcal{Z}$. It is not at all obvious why dimension reduction should be possible for the Toeplitz algebra, which has no local tensor product decomposition. The argument of \cite[Proposition~2.9]{WinZac:dimnuc} yields an approximation of $\mathcal{T}$ by noncommutative partitions of unity with three colours, hence showing that the nuclear dimension is at most two. (Of these three colours, two come from approximations of the circle; the third comes from the compacts.) In this paper we show that these approximations can be modified and repackaged so that only two colours suffice. One can think of our method as a noncommutative (and 2-coloured) bricklaying scheme; for this to work we will have to use the little available noncommutative space most efficiently. The key will be Lin's theorem on almost normal elements of matrix algebras; cf.\ \cite[The Main Theorem]{Lin:almost}.\\

The remainder of the paper is organised as follows: We first recall the definition of nuclear dimension from \cite{WinZac:dimnuc}. Then we apply the definition to the circle, both to make it explicit and to settle notation. Afterwards we isolate a lemma which encapsulates the additional geometric insight that makes everything work. Finally, we prove our result by exhibiting a system of approximations for the Toeplitz algebra that shows its nuclear dimension to be one.\\

\begin{dfn}
The nuclear dimension of a $\textup{C}^*$-algebra $A$ is the smallest natural number $d$ so that the following holds: 
There exists a net $(F_\lambda, \psi_\lambda, \phi_\lambda)_{\lambda \in \Lambda}$ such that the $F_\lambda$ are finite-dimensional $\textup{C}^*$-algebras, and such  that $\psi_\lambda \colon A \to F_\lambda$ and $\phi_\lambda \colon F_\lambda \to A$ are completely positive maps satisfying  
\begin{enumerate}[label=(\roman{*})]
\item  $\|\phi_\lambda \circ \psi_\lambda(a) - a\| \longrightarrow 0$ for all  $a \in A$, 
\item $\Vert \psi_\lambda \Vert \leq 1$,   
\item for each $\lambda$, $F_\lambda$ decomposes into $d+1$ ideals $F_\lambda = F_\lambda^{(0)} \oplus \ldots \oplus F_\lambda^{(d)}$ such that $\phi_\lambda^{(l)} := \phi_\lambda |_{F_\lambda^{(l)}}$ is a completely positive contractive order zero map for $l=0, ..., d$, where {\it order zero} means that $\phi^{(l)}_\lambda$ sends pairs of orthogonal elements to pairs of orthogonal elements. 
\end{enumerate}
\end{dfn}

\begin{ex}
The nuclear dimension of $A = C(S^1)$ is one, and can be realised  by completely positive approximations coming from `2-coloured' partitions of unity. We set these up explicitly, as we will need them later on in the paper.

For $0 \neq k \in \mathbb{N}$, $j = 1, \ldots, k$ and $l=0,1$ set $\delta_{j,k}^{(l)} := e^{2 \pi i \frac{2j+l}{2k}} \in S^1$. Define $\psi_k \colon C(S^1) \to \C^k \oplus \C^k$ by a sum of point evaluations, 
\[
\textstyle
\psi_k:= \psi_k^{(0)} \oplus \psi_k^{(1)} := \big(\bigoplus_{j=1}^k \textup{ev}_{\delta_{j,k}^{(0)}} \big)\oplus \big( \bigoplus_{j=1}^k \textup{ev}_{\delta_{j,k}^{(1)}} \big).
\] 

Again for $0 \neq k \in \mathbb{N}$, $j = 1, \ldots, k$ and $l=0,1$ take a partition of unity for $S^1$ by sawtooth functions $\{\bar{e}_{j,k}^{(l)}\}_{j \in \{1,\ldots k\}, \,l\in\{0,1\}}$ as depicted below:

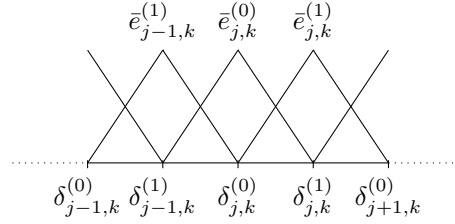
\[
\begin{tikzpicture}[domain=0:4]
% Achsen zeichnen
\draw[dotted] (2,0) -- (3,0) ;
\draw[-] (3,0) -- (7,0) ;
\draw[dotted] (7,0) -- (8,0) ;

\draw (3 cm,1pt) -- (3 cm,-2pt) node[anchor=north] {$\delta_{j-1,k}^{(0)}$};
\draw (4 cm,1pt) -- (4 cm,-2pt) node[anchor=north] {$\delta_{j-1,k}^{(1)}$};
\draw (4 cm,1.5) -- (4 cm,1.5) node[anchor=south] {$\bar{e}_{j-1,k}^{(1)}$};
\draw (5 cm,1pt) -- (5 cm,-2pt) node[anchor=north] {$\delta_{j,k}^{(0)}$};
\draw (5 cm,1.5) -- (5 cm,1.5) node[anchor=south] {$\bar{e}_{j,k}^{(0)}$};
\draw (6 cm,1pt) -- (6 cm,-2pt) node[anchor=north] {$\delta_{j,k}^{(1)}$};
\draw (6 cm,1.5) -- (6 cm,1.5) node[anchor=south] {$\bar{e}_{j,k}^{(1)}$};
\draw (7 cm,1pt) -- (7 cm,-2pt) node[anchor=north] {$\delta_{j+1,k}^{(0)}$};

\draw(3,1.5) -- (4,0);
\draw (3,0) -- (4, 1.5); 
\draw(4,1.5) -- (5,0);
\draw (4,0) -- (5, 1.5); 
\draw(5,1.5) -- (6,0);
\draw (5,0) -- (6, 1.5); 
\draw(6,1.5) -- (7,0);
\draw (6,0) -- (7, 1.5); 
\end{tikzpicture} 
\]

Define completely positive contractions $\phi_k^{(0)}, \phi_k^{(1)} \colon \C^k \to C(S^1)$ by
\begin{equation*}
\textstyle
\phi_k^{(l)}( a_1,  \ldots  ,a_k) :=  \sum_{j=1}^k a_j \bar{e}_{j,k}^{(l)} \text{ for } l=0,1,
\end{equation*}
and $\phi_k \colon \C^k \oplus \C^k \to C(S^1)$ by $\phi_k := \phi_k^{(0)} + \phi_k^{(1)}$; note that both maps  $\phi_k^{(l)}$  are order zero, i.e., preserve orthogonality. With this setup the sequence $(\phi_k \psi_k)_k$ converges to the identity map on $C(S^1)$ in point-norm topology.
\end{ex}

The lemma below is key for the proof of our main result. It relies on two nontrivial facts and a beautiful theorem. 

The first fact is \cite[Corollary~3.1]{WZ:MJM}: order zero contractions out of a C$^*$-algebra $A$ are in bijection with $^*$-homomorphisms out of $C_0((0,1]) \otimes A$, the cone over $A$. 

The second fact is that the cone over $C(S^1)$, which is isomorphic to $C_0(D \setminus \{0\})$ (with $D$ the unit disc), is canonically isomorphic to the universal C$^*$-algebra generated by a normal contraction. 

The theorem we will be using is Lin's result on almost normal matrices; it says that for every $\eta>0$ there is a $\delta>0$ such that, whenever $x$ is a square matrix with $\|x x^* - x^*x\| < \delta$, there is a square matrix $y$ such that $yy^* - y^*y = 0$ and $\|x - y\|< \eta$. The point is of course that $\delta$ depends only on $\eta$, and not on the matrix size. This was shown by Lin in \cite{Lin:almost}; an alternative proof was given by Friis and R{\o}rdam in \cite{FriisRor}. 

\begin{lem}
Let $(B_n)_{n \in \mathbb{N}}$ be a sequence of finite-dimensional  $\textup{C}^*$-algebras  and let    
\begin{equation*}
\textstyle
\beta \colon C(S^1) \longrightarrow \prod_{n \in \N} B_n \Big/ \sum_{n \in \N} B_n 
\end{equation*}
be a completely positive contractive order zero map. Then $\beta$ has a completely positive contractive order zero lift  
\begin{equation*}
\textstyle
\bar{\beta} \colon C(S^1) \longrightarrow \prod_{n \in \N} B_n.
\end{equation*}
\begin{proof}
Let $u \in C(S^1)$ be the identity function. By \cite[Corollary~3.1]{WZ:MJM}, $\beta$ induces a $^*$-homomorphism 
\begin{equation*}
\textstyle
\pi \colon  C_0((0, 1 ]) \otimes C(S^1) \longrightarrow  \prod_{n \in \N} B_n \Big/ \sum_{n \in \N} B_n 
\end{equation*}
with $\pi(\id_{(0,1]} \otimes u) = \beta(u)$. In particular we see that $\beta(u)$ is normal. Let $(v_n)_{n \in \mathbb{N}} \in \prod_{n \in \mathbb{N}} B_n$ be a contractive lift of $\beta(u)$. Note that each $v_n$ is a direct sum of matrices and that $\|v_n v_n^* - v_n^* v_n\|$ converges to $0$ as $i$ goes to infinity. But then it follows from Lin's theorem that there is a sequence $(\bar{v}_n)_{n \in \mathbb{N}} \in \prod_{n \in \mathbb{N}} B_n$ of {\it normal} contractions such that $\|\bar{v}_n - v_n\| \to 0$, i.e., $(\bar{v}_n)_{n \in \mathbb{N}}$ is a normal contractive lift of $\beta(u)$.

Upon identifying $C_0((0, 1 ]) \otimes C(S^1)$ with $C_0 (D \setminus \{ 0 \})$ and regarding the latter as the universal C$^*$-algebra generated by a normal contraction, we obtain a $^*$-homomorphism 
\begin{equation*}
\textstyle
\bar{\pi} \colon  C_0((0, 1 ]) \otimes C(S^1) \longrightarrow  \prod_{n \in \N} B_n
\end{equation*}
with $\bar{\pi}(\id_{(0,1]} \otimes u) = (\bar{v}_n)_{n \in \N}$. Note that this implies that $\bar{\pi}$ lifts $\pi$. Now define
\begin{equation*}
\textstyle
\bar{\beta}(\, .\,):= \bar{\pi}(\id_{(0,1]} \otimes \, . \,): C(S^1) \longrightarrow  \prod_{n \in \N} B_n,
\end{equation*}
then $\bar{\beta}$ is the composition of a $^*$-homomorphism with an order zero map, hence itself an order zero map. Moreover, since $\bar{\pi}$ lifts $\pi$, $\bar{\beta}$ lifts $\beta$.
\end{proof}
\end{lem}

\begin{thm}
The Toeplitz algebra has nuclear dimension one. 

\begin{proof}
First, we choose a quasicentral idempotent approximate unit $h= (h_n)_{n \in \N}$ for $\k \lhd \T$ so that every $h_n$ has finite rank (here, `idempotent' means that $h_{n+1} h_n = h_n$ for all $n \in \mathbb{N}$). Let $\tilde{h}=(\tilde{h}_n)_{n \in \N}$ be essentially the same approximate unit as $h$, but with index shifted by one (and with $h_0$ dropped). Upon regarding $h$ and $\tilde{h}$ as elements in $\prod_{n \in \N} \T$,  we have $h \cdot \tilde{h} = h$. 
\\ \\
To find approximations for $\mathcal{T}$, we work with the following diagram:

\begin{tikzpicture}[node distance = 2.5cm, auto]
  \node (T) {$\T$};
  \node (S) [right of=T] {$C(S^1)$};
  \node (TQ)  [node distance=9cm, right of=S]  {$\faktor{\prod \T} {\sum \T} $}; 
  \node(T2) [node distance=2cm, left of=TQ] {$\prod \T $} ; 
  \node(C1) [node distance= 2cm, above of= S]{$\C^k \oplus \C^k$};
  \node (C) [right of=C1, above of= S] {$\prod C_n$};
  \node (B) [node distance= 5cm, , below of=C] {$\prod B_n$};
  \node (A) [below of=B] {$\prod A_n$};
  \node (CQ) [right of=C]{$\faktor{\prod C_n} {\sum C_n} $}; 
  \node (BQ) [right of=B]{$\faktor{\prod B_n} {\sum B_n} $}; 
  \node (C2) [right of=S]{$\C^k$};
  \draw[->] ([yshift=2pt]  T.0) to node {$\pi$} ([yshift=2pt]  S.180)   ;
  \draw[->] ([yshift=-2pt] S.180) -- node {$\mu$} ([yshift=-2pt] T.0);  
  \draw[->] (C) to node {$q_C $} (CQ) ; 
  \draw[->] (B) to node {$q_B $} (BQ) ;
  \draw[->] (T2) to node {$q$} (TQ);
  \draw[->, bend right] (T) to node {$\alpha $} (A) ;
  \draw[->] (C) to node {$\iota_4$} (T2);
  \draw[->] (CQ) to node {$\iota_3$} (TQ);
  \draw[->] (S) to node {$\gamma $} (C);
  \draw[->, dashed] ([xshift=-4pt] S.90) to node {$\psi_k $} ([xshift=-4pt] C1.270);
  \draw[->, dashed] ([xshift=4pt] C1.270) to node {$\phi_k $} ([xshift=4pt] S.90) ;
  \draw[->, dashed] (S) to node {$\psi_k^{(0)} $} (C2) ;
  \draw[->, dashed] (C2) to node {$\rho_k $} (B);
  \draw[->, thick] (S) to node {$\bar{\beta} $} (B);
  \draw[->, bend right] (S) to node[swap] {$\beta $} (B);
  \draw[->, bend left=70] (T) to node {$\iota $} (TQ);
  \draw[->, name path = p1 , bend right= 42] (A) to node[swap] {$\iota_1$} (T2);
  \path[name path=p2] (BQ) --(TQ);	
  %add white space
   \path [name intersections={of = p1 and p2}];
  \coordinate (S)  at (intersection-1);
  \path[name path=circle] (S) circle(1mm);
  % find intersections of second line and circle
  \path [name intersections={of = circle and p2}];
  \coordinate (I1)  at (intersection-2);
  \coordinate (I2)  at (intersection-1);
   \draw (BQ) to node  {$\iota_2$} (I1);
   \draw[->] (I2) -- (TQ);
  \end{tikzpicture}

Here $\pi$ is the quotient map from the short exact sequence, and $\mu$ is a completely positive contractive lift for $\pi$ (for example map $u \in C(S^1)$ to the bilateral shift on $\ell^2(\mathbb{Z})$ and compress to $\ell^2(\mathbb{N})$). Moreover, we let 
\begin{equation*}
A_n = h_n \T h_n, \quad B_n = \big(\tilde{h}_n - h_n \big) \T \big(\tilde{h}_n - h_n \big),  \quad C_n = \overline{\big(1-\tilde{h}_n\big) \T \big(1-\tilde{h}_n \big)} 
\end{equation*}
for all $n \in \N$, and 
\begin{align*}
&\alpha \colon x \longmapsto \big(h_n^{\frac{1}{2}} x h_n^{\frac{1}{2}}\big)_{n \in \N}, \\
&\beta \colon f \longmapsto \big((\tilde{h}_n - h_n)^{\frac{1}{2}} \mu(f) (\tilde{h}_n - h_n)^{\frac{1}{2}} \big)_{n \in \N},  \\
&\gamma \colon f \longmapsto \big((1-\tilde{h}_n)^{\frac{1}{2}} \mu(f) (1-\tilde{h}_n)^{\frac{1}{2}} \big)_{n \in \N}. 
\end{align*}
Let $q, q_B, q_C$ be the canonical quotient maps, let $\iota_1, \iota_2, \iota_3, \iota_4$ be the natural inclusion maps, and let $\iota$ be the canonical constant sequence  embedding. 

The map $\bar{\beta}$ is a completely positive contractive order zero lift of the order zero map $q_B \circ \beta$; this is where our lemma above (and hence Lin's theorem) enters.  

Before explaining the dotted arrows, note that with the setup so far the map $\iota$ is precisely the sum of the lower three paths from left to right through the diagram,
\begin{equation}
\label{iota-sum}
\iota = q \circ \iota_1 \circ \alpha + \iota_2 \circ q_B \circ \bar{\beta} \circ \pi + \iota_3 \circ q_C \circ \gamma \circ \pi. 
\end{equation}

For each $0 \neq k \in \mathbb{N}$, let  
\begin{equation*}
C(S^1) \stackrel{\psi_k}{\longrightarrow} \mathbb{C}^k \oplus \mathbb{C}^k \stackrel{\phi_k}{\longrightarrow} C(S^1)
\end{equation*}
with $\psi_k = \psi_k^{(0)} \oplus \psi_k^{(1)}$ and $\phi_k = \phi_k^{(0)} + \phi_k^{(1)}$
be the completely positive approximations described in our example above.

The definition of the maps $\rho_k$ in the diagram above is more complicated, but it is at the heart of the matter. We define them component-wise, working with the diagram below: 
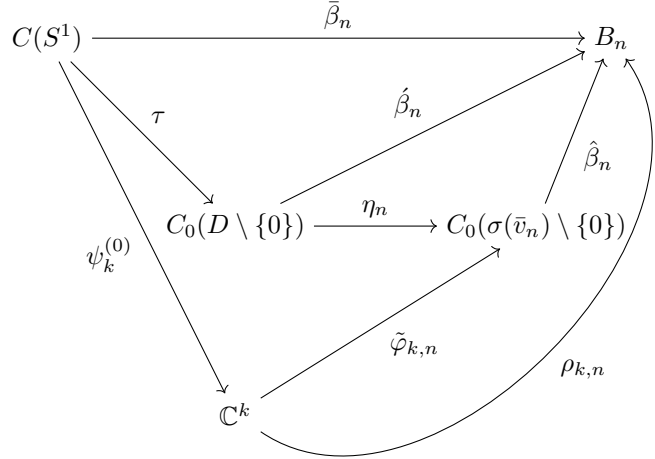
\[
\begin{tikzpicture}[node distance=2.5cm, auto]
\node (S)  {$C(S^1)$};
\node (B) [node distance= 7.5cm, right of=S] {$ B_n$};
\node (D) [right of= S, below of=S]{$C_0(D \setminus \{ 0 \})$};
\node (C1)[below of=D] {$\C^k$} ; 
\node (v) [node distance= 4cm, right of=D]{$C_0(\sigma(\bar{v}_n) \setminus \{ 0 \} )$ }; 
\draw[->] (S) to node {$\bar{\beta}_n $} (B);
\draw[->] (S) to node {$\tau $} (D);
\draw[->] (S) to node[swap] {$\psi_k^{(0)} $} (C1) ;
\draw[->] (D) to node {$\eta_n $} (v);
\draw[->] (v) to node[swap] {$\hat{\beta}_n$} (B) ;
\draw[->] (D) to node {$\acute{\beta}_n $} (B) ;
\draw[->] (C1) to node[swap] {$\tilde{\phi}_{k, n} $} (v) ;
\draw[->, bend right=80] (C1) to node[swap] {$\rho_{k,n}$} (B);
\end{tikzpicture}
\]
Here, $\tau$ is the canonical order zero map which sends the identity function $u \in C(S^1)$ to the identity function on $D$ (this is just the map $(f \mapsto \id_{(0,1]} \otimes f)$ followed by the identification $C_0((0,1]) \otimes C(S^1) \cong C_0(D \setminus \{0\})$ given by $(\mathrm{id}_{(0,1]} \otimes u \longmapsto \mathrm{id}_D)$. 

We write $\bar{\beta}_n$ for the (order zero) components of $\bar{\beta}$. With the notation of our lemma, we have normal elements $\bar{\beta}_n(u) = \bar{v}_n$, $n \in \mathbb{N}$. Upon identifying $C_0(D \setminus \{0\})$ with the universal C$^*$-algebra generated by a normal contraction, we obtain the canonically induced $^*$-homomorphisms $\acute{\beta}_n$, which canonically factorise through continuous functions on the spectrum of $\bar{v}_n$ via the $^*$-homomorphisms $\eta_n$ and $\hat{\beta}_n$.

For the definition of the maps $\tilde{\phi}_{k, n}$ and $\rho_{k,n}$ we write the pointed disc $D \setminus \{0\}$ as the disjoint union of `pizza slices' 
\begin{equation*}
\textstyle
S_{j,k}:= \big\{r e^{2\pi i t} \mid \frac{ 2j-1}{2k} < t \le \frac{ 2j +1}{2k} ; \; 0<r\le 1 \big\} \mbox{ for } j=1,\ldots,k.
\end{equation*}
Let $\chi_{S_{j,k}}$ be the characteristic functions of these sets and define Borel functions  $\acute{\chi}_{j,k}: D  \longrightarrow \mathbb{C}$ by
\begin{equation*}
\acute{\chi}_{j,k}(z):= \chi_{S_{j,k}}(z) \cdot |z|.
\end{equation*}
Now since for each $n$ the spectrum of the element $\bar{v}_n \in B_n$ is discrete, the $\acute{\chi}_{j,k}$ restrict to \emph{continuous} functions on $\sigma(\bar{v}_n)$ and we may define completely positive contractive order zero maps $\tilde{\phi}_{k,n}:\mathbb{C}^k \longrightarrow C_0(\sigma(\bar{v}_n) \setminus \{0\})$ by 
 \begin{equation*}
 \textstyle
\tilde{\phi}_{k, n} \colon (a_1, \ldots, a_k) \longmapsto   \sum_{j=1}^{k} a_j  \cdot \acute{\chi}_{j,k} |_{\sigma(\bar{v}_n) \setminus \{ 0 \} }.
\end{equation*}
We then define 
\begin{equation*}
\rho_{k,n}:= \hat{\beta}_n \circ \tilde{\phi}_{k,n}
\end{equation*} 
and compute
\begin{align*}
\textstyle
\rho_{k,n} \circ \psi_k^{(0)} (f) & \textstyle = \hat{\beta}_n \big( \sum_{j=1}^k  f(\delta^{(0)}_{j,k}) \cdot  \acute{\chi}_{j,k}|_{\sigma(\bar{v}_n) \setminus\{0\}} \big) \\
& \textstyle = \hat{\beta}_n \circ \eta_n \big( \sum_{j=1}^k  f(\delta^{(0)}_{j,k}) \cdot  \acute{\chi}_{j,k}  \big) 
\end{align*}
for $f \in C(S^1)$. Here, the $\delta^{(0)}_{j,k}$ are as in our example above, and we have we have tacitly extended the restriction $^*$-homomorphism $\eta_n$ to all bounded Borel functions on $D$.

Now observe that for each $f \in C(S^1)$ we have
\begin{equation*}
\textstyle
\big\|\sum_{j=1}^k  f(\delta^{(0)}_{j,k}) \cdot  \acute{\chi}_{j,k}  - \tau(f) \big\|_{\infty,D} \stackrel{k \to \infty}{\longrightarrow} 0,
\end{equation*}
whence 
\begin{align*}
\big\| \rho_{k,n} \circ \psi_k^{(0)}(f) - \bar{\beta}_n(f) \big\| & \textstyle = \big\| \hat{\beta}_n \circ \eta_n  \big( \sum_{j=1}^k  f(\delta^{(0)}_{j,k}) \cdot  \acute{\chi}_{j,k}  \big)  -  \hat{\beta}_n \circ \eta_n \circ \tau(f) \big\|  \\
& \textstyle \le \big\| \big( \sum_{j=1}^k  f(\delta^{(0)}_{j,k}) \cdot  \acute{\chi}_{j,k} \big) - \tau(f) \big\|_{\infty,D} \\
& \stackrel{k \to \infty}{\longrightarrow} 0
\end{align*}
for each $f \in C(S^1)$. Note that convergence here does not depend on $n$. As a consequence, we see that the maps $\rho_k \circ \psi^{(0)}_k$ converge to $\bar{\beta}$ in point-norm topology,
\begin{equation*}
\big\|\rho_k \circ \psi^{(0)}_k (f) - \bar{\beta} (f) \big\| \stackrel{k \to \infty}{\longrightarrow} 0
\end{equation*}
for every $f \in C(S^1)$. Upon plugging this into \eqref{iota-sum} we obtain for every $x \in \mathcal{T}$
\begin{equation}
\label{sum2}
\big\Vert \iota(x) - q  \, \iota_1 \,  \alpha(x) + \iota_2 \, q_B \, \rho_k \, \psi_k^{(0)} \, \pi(x) + \iota_3 \, q_C \, \gamma \, \phi_k \, \psi_k \, \pi(x) \big\Vert \stackrel{k \to \infty}{\longrightarrow} 0.
\end{equation}

Next, we check that for each $k$ the completely positive contractive map
\begin{equation*}
\textstyle
\iota_2 \circ q_B \circ \rho_k + \iota_3 \circ q_C \circ \gamma \circ \phi_k^{(0)} \colon \C^k \longrightarrow  \prod \T \big/ \sum \T 
\end{equation*}
is order zero. Since the two summands are each order zero, we only need to examine the mixed terms, i.e., we have to confirm that
\begin{equation}
\label{orthogonality}
\iota_2 \circ q_B \circ \rho_k (e_{j})  \cdot  \iota_3 \circ q_C \circ \gamma \circ \phi_k^{(0)}(e_{j'})= 0 \text{ for } 1 \le j \neq j' \le k,
 \end{equation} 
where we write $e_j$ for the canonical $j$th generator of $\mathbb{C}^k$. The reason why \eqref{sum2} holds is that our `pizza slices' $\acute{\chi}_{j,k}$ and the $\tau (\bar{e}^{(0)}_{j',k})$ are orthogonal functions over $D$ whenever $j \neq j'$. To be more precise, let $(d_m)_{m \in \mathbb{N}} \subset C(S^1)$ be normalised functions such that $d_m(e^{2 \pi i t}) =1 $ for all $\frac{2j-1}{2k} < t \le \frac{2j+1}{2k}$ and such that $\bar{e}^{(0)}_{j',k} \cdot d_m$ goes to zero uniformly as $m$ goes to infinity. We then have $\acute{\chi}_{j,k} \le \tau(d_m)$, whence 
\begin{equation*}
\rho_{k,n}(e_j) = \hat{\beta}_n \circ \tilde{\phi}_{k,n}(e_j) = \hat{\beta}_n (\acute{\chi}_{j,k}|_{\sigma(\bar{v}_n) \setminus \{0\}})\le \hat{\beta}_n \circ \eta_n \circ \tau (d_m) = \bar{\beta}_n(d_m)
\end{equation*}
for each $m$ and $n$. It follows that 
\begin{equation*}
q_B \circ \rho_k(e_j) \le q_B \circ \bar{\beta} (d_m) = q_B \circ \beta (d_m)
\end{equation*}
for each $m$, and therefore
\begin{align*}
\lefteqn{\big\| \iota_3 \circ q_C \circ \gamma \circ \phi_k^{(0)}(e_{j'}) \cdot \iota_2 \circ q_B \circ \rho_k (e_{j}) \big\|^2} \\
& \le \big\| \iota_3 \circ q_C \circ \gamma \circ \phi_k^{(0)}(e_{j'}) \cdot \iota_2 \circ q_B \circ \beta(d_m) \cdot \iota_3 \circ q_C \circ \gamma \circ \phi_k^{(0)}(e_{j'}) \big\| \\
& \le \big\| (1 - \tilde{h})^{\frac{1}{2}} \iota(\mu(\bar{e}^{(0)}_{j',k}))(1 - \tilde{h})^{\frac{1}{2}}  \cdot (\tilde{h} - h)^{\frac{1}{2}} \iota(\mu(d_m)) (\tilde{h} - h)^{\frac{1}{2}} \big\| \\
& = \big\| (1 - \tilde{h}) (\tilde{h} - h)\iota(\mu(\bar{e}^{(0)}_{j',k})) \cdot \iota(\mu(d_m))  \big\| \\
& \le \big\| \iota(\mu(\bar{e}^{(0)}_{j',k} \cdot d_m )) \big\| \\
& \stackrel{m \to \infty}{\longrightarrow} 0 ,
\end{align*} 
which entails \eqref{orthogonality}.

Next, recall that by \cite[Remark~2.4]{KirchbergWinter:CovDim} order zero maps out of finite-dimensional C$^*$-algebras into quotient C$^*$-algebras always lift to order zero maps, and so there are completely positive contractive order zero lifts $\tilde{\phi}_{k}^{(1)} \colon \C^k \longrightarrow \prod C_n$ of $q_C \circ \gamma \circ \phi_k^{(1)}$, and $\tilde{\phi}_k^{(0)} \colon \C^k \longrightarrow \prod \T$ of $\iota_2 \circ q_B \circ \rho_k + \iota_3 \circ q_C \circ \gamma \circ \phi_k^{(0)}$, fitting into the following diagram: 
\[
\begin{tikzpicture}[node distance=2.5cm, auto]
\node (T) {$\T$};
\node (S) [right of=T] {$C(S^1)$};
\node (C1)[node distance= 1.5cm, right of=S, above of=S] {$\C^k$} ; 
\node (C2)[node distance= 1.5cm, right of=S, below of=S] {$\C^k$} ; 
\node (T1) [node distance = 6.5cm, right of=S ] {$\prod \T$};
\node (TQ)  [node distance=8.5cm, right of=S]  {$\faktor{\prod \T} {\sum \T} $}; 
\node (A) [below of=C2] {$\prod A_n$};
\node (C) [  right of=C1, above of= C2] {$\prod C_n$};
\draw[->] (T) to node {$\pi $} (S);
\draw[->] (S) to node {$\psi_k^{(1)} $} (C1) ;
 \draw[->] (S) to node[swap] {$\psi_k^{(0)}$} (C2);
 \draw[->] (T) to node[swap] {$\alpha$} (A) ;
 \draw[->] (A) to node[swap] {$ \iota_1$} (T1);
 \draw[-> ] (C1) to node {$\tilde{\phi}_{k}^{(1)}$} (C);
 \draw[->] (C) to node {$\iota_4 $} (T1);
 \draw[->] (C2) to node {$\tilde{\phi}_k^{(0)} $} (T1);
 \draw[->] (T1) to node {$q$} (TQ);
 \draw[->, bend left=50] (T) to node {$\iota $} (TQ);
\end{tikzpicture}
\]
We already know from \eqref{sum2} that for every $x \in \mathcal{T}$
\begin{equation*}
\big\| \iota(x) - q  \big(\iota_1 \circ \alpha(x) +  \tilde{\phi}_k^{(0)} \circ \psi_k^{(0)} \circ \pi(x) +  \iota_4 \circ \tilde{\phi}^{(1)}_k \circ \psi^{(1)}_k \circ \pi (x)\big)\big\| \stackrel{k \to \infty}{\longrightarrow} 0;
\end{equation*}
expressing the quotient norm as a limit we obtain for every $x \in \mathcal{T}$
\begin{equation*}
\lim_{k \to \infty} \lim_{n \to \infty} \big\| x - \iota_{1,n} \circ \alpha_n(x) +  \tilde{\phi}_{k,n}^{(0)} \circ \psi_{k}^{(0)} \circ \pi(x) +  \iota_{4,n} \circ \tilde{\phi}^{(1)}_{k,n} \circ \psi^{(1)}_k \circ \pi (x) \big\| = 0.
\end{equation*}
This remains true for finite sets of elements $x$, and using separability of $\mathcal{T}$ we run a diagonal sequence argument to find, for every index $k$, some $n_k$ such that for $x \in \mathcal{T}$
\begin{equation*}
\lim_{k \to \infty} \big\| x - \iota_{1,n_k} \circ \alpha_{n_k}(x) +  \tilde{\phi}_{k,n_k}^{(0)} \circ \psi_{k}^{(0)} \circ \pi(x) +  \iota_{4,n_k} \circ \tilde{\phi}^{(1)}_{k,n_k} \circ \psi^{(1)}_k \circ \pi (x)) \big\| = 0.
\end{equation*}
Now for each $0\neq k \in \mathbb{N}$ we set 
\begin{align*}
&\dot{F}^{(0)}_k := \C^k , \\
&\dot{F}^{(1)}_k :=  \C^k \oplus A_{n_k}, \\
&\dot{F}_k := \dot{F}_k^{(0)} \oplus \dot{F}_k^{(1)} = \C^k \oplus \C^k \oplus A_{n_k}, \\
&\dot{\psi}_k := \psi_k^{(0)} \circ \pi + \psi_k^{(1)} \circ \pi + \alpha_{n_k} \colon  \T \longrightarrow \C^k \oplus \C^k \oplus A_{n_k}, \\
&\dot{\phi}_k^{(0)}  := \tilde{\phi}_{n_k}^{(0)} \colon \C^k \longrightarrow \T, \\
&\dot{\phi}_k^{(1)} := \iota_{4,n_k} \circ \phi^{(1)}_{C,n_k} \oplus \iota_{1,n_k} \colon \C^k \oplus A_{n_k} \longrightarrow \T, \\
&\dot{\phi}_k :=  \dot{\phi}_k^{(0)} \oplus \dot{\phi}_k^{(1)} \colon \C^k \oplus \C^k \oplus A_{n_k} \longrightarrow \T ;
\end{align*} 
note that $\phi_k^{(1)}$ is a completely positive order zero contraction since each of the two summands is, and since they map into the orthogonal algebras $C_{n_k}$ and $A_{n_k}$.

The system $(\dot{F}_k, \dot{\psi}_k, \dot{\phi}_k)_{k \in \mathbb{N} \setminus \{0\}}$ indeed establishes that the Toeplitz algebra has nuclear dimension one.
\end{proof}
\end{thm}

\bibliographystyle{amsplain}
\bibliography{Toeplitz-dimnuc}
 
\bigskip

\end{document}